\documentclass[11pt, a4paper]{article}
\usepackage{times}
\usepackage{a4wide}
\usepackage{hyperref}
\usepackage[british]{babel}
\usepackage{enumerate, longtable}
\usepackage{amsmath, amscd, amsfonts, amsthm, amssymb, latexsym, comment, graphicx, color, url}
\usepackage[all]{xypic}
\usepackage[T1]{fontenc}
\usepackage[utf8]{inputenc}
\usepackage{verbatim}
\usepackage{blindtext}
\usepackage{multicol}
\usepackage{multicap}

\newtheorem{thm}{Theorem}[section]

\newtheorem{defi}[thm]{Definition}

\newcommand{\donothing}[1]{}

\newcommand{\QQ}{\mathbb{Q}}

\newcommand{\ZZ}{\mathbb{Z}}

\newcommand{\exclude}[1]{}

\begin{document}

\selectlanguage{british}

\title{On The Hecke Eigenforms of Half-Integral Weight and Dedekind-Eta Products}
\author{Banu Irez Aydin \footnote{Bilecik Seyh Edebali University, Vocational School, 11200 Bilecik, Turkey, banu.irez@bilecik.edu.tr}  and Ilker Inam\footnote{Bilecik Seyh Edebali University, Department of Mathematics, Faculty of Arts and Sciences, 11200 Bilecik, Turkey, ilker.inam@gmail.com or ilker.inam@bilecik.edu.tr }}
\maketitle

\begin{abstract}
Systematic choice of the Hecke eigenforms of half-integral weight is an interesting problem in the theory of modular forms. In this paper, we find all Dedekind-eta products of half-integral weight which are Hecke eigenforms up to weight $15/2$ with varying levels. Proof is based on the Shimura lift.

\end{abstract}

\section{Introduction}

In one hand, with various applications in different branches of mathematics and even physics modular forms have been receiving much attention for decades. Since space of modular forms a finite dimensional vector space over complex field, it is natural to look for eigenvectors under some linear maps, here so-called Hecke operators, therefore we have an important class of modular forms: Hecke eigenforms. They are really important. Indeed, there are many reasons for that. For instance, the coefficient function $a_f(n)$ is a multiplicative arithmetic function
they form a natural basis (of the newspace)
they have $L$-functions with Euler products, analytic/meromorphic continuation if weight 2 and coefficients in $\QQ$, those $L$-functions match with those of elliptic curves. Furthermore, they have attached Galois representations. Their Fourier coefficients are multiplicative and satisfy simple recursions at powers of primes. One of the breakthrough result of the century: Sato-Tate Conjecture is valid for Hecke eigenforms. Modular forms are computational-friendly objects and the reader is refer to \cite{BeCo} for both theoric and computational aspects of modular forms.

On the other hand Dedekind eta quotients are also very important. They are very simple examples for modular forms. Even in physics, for example, the (bosonic) string partition function can be expressed in terms of Dedekind eta function. As a historical note, in 1742, Euler proposed the following amazing identity:

$$ \prod_{n=1}^{\infty} (1-q^n) = \sum_{m=-\infty}^{\infty} (-1)^m  q^{\frac{1}{2} m(3m-1)} $$

And he succeeded to prove it later in 1750. This leads to the appearance of the Dedekind-eta function. There are many examples for applications of the Dedekind-eta quotients.

Systematic choice of Hecke eigenforms of half-integral weight is an interesting problem and there are some recent results in \cite{IW-Fast} and \cite{TJM} where Rankin-Cohen brackets of Eisenstein series and theta series are used effectively. The problem of fast computation of Fourier coefficients of $4$ level and half-integral weight Hecke eigenforms has been addressed in \cite{IW-Fast} and the systematic selection of $4N$-level and half-integral weight Hecke eigenforms with $N>1$ is still an open problem in the literature. In another direction, here, it is aimed to give a partial classification of Hecke eigenforms in terms of Dedekind-eta quotients. 

In this paper, Hecke eigenforms of half integral weight for all possible weights up to $15/2$ and various levels are expressed in terms of Dedekind eta products. A computerized proof based on the Shimura lift is employed and we give all necessary codes for the reader. 

\section*{Acknowledgment} The authors would like to thank Henri Cohen for the project in Luxembourg Summer School in 2018 and all interesting discussions. We are also grateful to John Voight for the kind help.

\section{Main Result}

Here is the main result of the paper.

\begin{thm}\label{main-thm}
All eta products listed in Table \ref{Table:1} are Hecke eigenforms of indicated weight, level and character. In other words, they are all eta products which are Hecke eigenforms of half-integral weight.

\newpage

\begin{table}[ht]
\centering
\begin{multicols}{2}
\noindent  
\scalebox{0.63}{

\begin{tabular}{|l||r|r|r|r|}
\hline

\hline
Example & Level & Weight & Character & Eta Quotient \\
\hline
1) & 576 & 1/2 & 12 & $\eta(24z)$ \\
\hline
2) & 1152 & 1/2 & 24 & $\eta(48z)$ \\ 
\hline
3) & 176 & 3/2 & 44 & $\eta^2(z)\eta(22z)$ \\
\hline
4) & 160 & 3/2 & 40 & $\eta^2(2z)\eta(20z)$ \\
\hline
5) & 864 & 3/2 & 8 & $\eta(2z)\eta(4z)\eta(18z)$ \\
\hline
6) & 448 & 3/2 & 28 & $\eta(2z)\eta(8z)\eta(14z)$ \\
\hline
7) & 176 & 3/2 & 1 & $\eta(2z)\eta^2(11z)$ \\
\hline
8) & 432 & 3/2 & 1 & $\eta^2(3z)\eta(18z)$ \\
\hline
9) & 864 & 3/2 & 8 & $\eta(3z)\eta(9z)\eta(12z)$ \\
\hline
10) & 128 & 3/2 & 8 & $\eta^2(4z)\eta(16z)$ \\
\hline
11) & 128 & 3/2 & 8 & $\eta^3(16z)$ \\
\hline
12) & 576 & 3/2 & 12 & $\eta(4z)\eta(8z)\eta(12z)$ \\
\hline
13) & 160 & 3/2 & 8 & $\eta(4z)\eta^2(10z)$ \\
\hline
14) & 288 & 3/2 & 24 & $\eta^2(6z)\eta(12z)$ \\
\hline
15) & 432 & 3/2 & 12 & $\eta(6z)\eta^2(9z)$ \\
\hline
16) & 64 & 3/2 & 1 & $\eta^3(8z)$ \\
\hline
17) & 720 & 5/2 & 5 & $\eta^3(z)\eta(6z)\eta(15z)$ \\
\hline
18) & 28 & 5/2 & 28 & $\eta^2(z)\eta^2(4z)\eta(14z)$ \\
\hline
19) & 128 & 5/2 & 8 & $\eta^4(2z)\eta(16z)$ \\
\hline
20) & 224 & 5/2 & 56 & $\eta^3(2z)\eta(4z)\eta(14z)$ \\
\hline
21) & 96 & 5/2 & 8 & $\eta^3(2z)\eta(6z)\eta(12z)$ \\
\hline
22) & 288 & 5/2 & 24 & $\eta^2(2z)\eta^2(4z)\eta(12z)$ \\
\hline
23) & 16 & 5/2 & 8 & $\eta^2(2z)\eta(4z)\eta^2(8z)$ \\
\hline
24) & 96 & 5/2 & 24 & $\eta(2z)\eta(4z)\eta^3(6z)$ \\
\hline
25) & 288 & 5/2 & 24 & $\eta^4(3z)\eta(12z)$ \\
\hline

\end{tabular}

\noindent 
\begin{tabular}{|l||r|r|r|r|}
\hline
Example & Level & Weight & Character & Eta Quotient \\
\hline
26) & 432 & 5/2 & 1 & $\eta^3(3z)\eta(6z)\eta(9z)$ \\
\hline
27) & 144 & 5/2 & 12 & $\eta^2(3z)\eta^3(6z)$ \\
\hline
28) & 64 & 5/2 & 1 & $\eta^4(4z)\eta(8z)$ \\
\hline
29) & 288 & 5/2 & 8 & $\eta^3(4z)\eta^2(6z)$ \\
\hline
30) & 144 & 7/2 & 1 & $\eta^3(z)\eta(3z)\eta^3(6z)$ \\
\hline
31) & 48 & 7/2 & 12 & $\eta^2(z)\eta^2(2z)\eta^3(6z)$ \\
\hline
32) & 64 & 7/2 & 1 & $\eta^4(2z)\eta^2(4z)\eta(8z)$ \\
\hline
33) & 48 & 7/2 & 1 & $\eta^3(2z)\eta^2(3z)\eta^2(6z)$ \\
\hline
34) & 32 & 7/2 & 8 & $\eta^2(z)\eta^5(4z)$ \\
\hline
35) & 144 & 7/2 & 12 & $\eta^3(4z)\eta^2(6z)$ \\
\hline
36) & 48 & 9/2 & 12 & $\eta^6(z)\eta^3(6z)$ \\
\hline
37) & 32 & 9/2 & 8 & $\eta^4(z)\eta^5(4z)$ \\
\hline
38) & 80 & 9/2 & 5 & $\eta^3(z)\eta^3(2z)\eta^3(5z)$ \\
\hline
39) & 16 & 9/2 & 1 & $\eta^2(z)\eta^3(2z)\eta^4(4z)$ \\
\hline
40) & 64 & 9/2 & 1 & $\eta^8(z)\eta(8z)$ \\
\hline
41) & 32 & 9/2 & 8 & $\eta^6(2z)\eta^3(4z)$ \\
\hline
42) & 48 & 9/2 & 1 & $\eta^3(2z)\eta^6(3z)$ \\
\hline
43) & 4 & 11/2 & 1 & $\eta^2(z)\eta^7(2z)\eta^2(4z)$ \\
\hline
44) & 32 & 11/2 & 8 & $\eta^{10}(2z)\eta(4z)$ \\
\hline
45) & 48 & 13/2 & 12 & $\eta^5(z)\eta^5(2z)\eta^3(3z)$ \\
\hline
46) & 16 & 13/2 & 1 & $\eta^2(z)\eta^11(2z)$ \\
\hline
47) & 16 & 15/2 & 1 & $\eta^6(z)\eta^9(2z)$ \\
\hline
48) & 48 & 13/2 & 12 & $\eta^5(z)\eta^5(2z)\eta^3(3z)$ \\
\hline
49) & 16 & 13/2 & 1 & $\eta^2(z)\eta^11(2z)$ \\
\hline
50) & 16 & 15/2 & 1 & $\eta^6(z)\eta^9(2z)$ \\
\hline
\end{tabular}

}

\end{multicols}
\caption{Hecke eigenforms in terms of Dedekind-eta quotients}\label{Table:1}
\end{table}

\end{thm}

\section{Preliminaries and Background Material}

In this section, we will give the necessary background. Let us begin with the definition of half-integral weight modular forms.

\begin{defi}
A modular form of half-integral weight $k/2$ , level $N$ and character $\chi$ is an holomorphic function defined on the upper half plane $\mathcal{H}$, satisfying the transformation formula
$$
f(\gamma z) = \chi(d) j(\gamma, z)^k f(z)
$$
for every $\gamma=\left(\begin{matrix} a & b \\ c & d \end{matrix}\right)\in\Gamma_0(4N)$ and $z\in \mathcal{H}$, and being holomorphic at the cusps. Here the automorphy factor $j(\gamma, z)$ is given by
$$
j(\gamma, z) = \theta(\gamma z)/\theta(z),
$$
where $\theta$ denotes the classical theta function. This space is denoted by $M_{k/2}(N)$, and the subspace of cusp forms is denoted by $S_{k/2}(N)$.
\end{defi}

Here is the definition of the Dedekind-eta function.

\begin{defi}
It is well-known that the Dedekind-eta function is defined by 
$$\eta(z)=q^{1/24} \prod_{n\geq 1} (1-q^n) $$

with $q=e^{2 \pi iz}$ and $Im(z)>0$. 

It is well-known that the Dedekind-eta function is a modular form of weight $1/2$ on $SL(2, \ZZ )$ with a complicated multiplier system formed of $24$th roots of unity. Its 24th power is the discriminant function $\Delta(z )$. 
\end{defi}

Here is the definition of Dedekind-eta product.
\begin{defi}
An eta product is any function of the form $$F(z) = \prod_{m \in I} {\eta (m z)^{r_m}}$$
where $I$ is a finite set of positive integers, and $r_m >0$ for all $m \in I$ and $Im(z)>0$.
\end{defi}

It is immediate to see that $\eta (m z)$ is modular on $\Gamma_0 (m)$ (with multiplier system), hence $F$ is modular on $\Gamma_0(M)$, where $M$ is the lowest common multiple of $m \in I$.\\

To obtain Hecke eigenforms, first we should define Hecke operators since Hecke eigenforms are exactly eigenvectors for all Hecke operators. The following theorem gives an explicit description of Hecke operators in half-integral weight setting. Note that, we can only define prime-squared-th Hecke operators for half-integral modular forms.

\begin{thm} \cite{Ko} Suppose that $4 \mid N$, $\chi$ is a Dirichlet character modulo $N$, $p \nmid N$ is a prime and $k=2\lambda+1$ is a positive odd integer. Let $f(z)=\sum_{n=0}^\infty  a_n q^n \in  M_{k/2}(\Gamma_0(N),\chi)$. Then $p^2$-th Hecke operator for $f$ is defined as
$$T_{p^2}f(z)=\sum_{n=0}^\infty b_n e^{2\pi inz}$$
where
$$b_n=a_{p^2n}+\chi(p)( \frac{(-1)^\lambda n}{p}) p^{\lambda-1}a_n + \chi(p^2)p^{k-2} a_{n/p^2}.$$

Here we take $a_{n/p^2}=0$ if $p^2 \nmid n$.

\end{thm}

Modular forms have important arithmetic properties and the most crucial step is of course to determine when two modular forms are equal. The Sturm bound fills the gap here. This bound is valid at both integral and half-integral weights. 

\begin{thm}
The Sturm bound for a space of modular forms $M_k(G,\chi)$ is a number $s$ such that if $f=\sum_{n\geqslant0} a_n q^n$ and $g=\sum_{n\geqslant0} b_n q^n$ are both elements in $M_k(G,\chi)$ and $a_n=b_n$ for $0\leqslant n \leqslant s+1$ then $f=g$.
\end{thm}

The following theorem gives necessary and sufficient condition for an eta product to become a modular (cusp) form. 

\begin{thm}\cite{CoSt}
Let $F=\prod_{m \in I} \eta(mz)^{r_m}$ be an eta product, so that $F$ has weight $k=\sum_{m \in I} r_m/2$, and denote by $M$ the least common multiple of the elements of $I$. Then $F$ belongs to some $M_k(\Gamma_{0}(N), \chi)$ if and only if $\sum_{m \in I} m r_m \equiv 0 (\mod 24)$, and if this is the case $F$ will be a cusp form. One can choose for $N$ the least common multiple of $M$ and the denominator of $\sum_{m \in I} r_m / (24m) $, and $\chi=\chi_D$ (quadratic character of discriminant $D$) with $D=(-1)^k \prod_{m \in I} m^{r_m}$ if $k \in \ZZ$, or $D=8 \prod_{m \in I} m^{r_m}$ if $k \in 1/2+ \ZZ$. 
\end{thm}

Main ingredient of proof of the main result will be the following result which is based on the Shimura lift.

\begin{thm}\cite{Purkait} \label{Soma}
Let $k, N$ be positive integers with $k \geq 3$ odd, $4 \mid N$. Let $\chi$ be a Dirichlet character modulo $N$. Let $N'=N/2$. Define

$$ m=N'^2 \prod_{p \mid N'} \left( 1- \frac{1}{p^2}  \right), R=\frac{(k-1)m}{12}-\frac{m-1}{N'} .$$

Then $T_{i^2}$ for $i \leq R$ generate the restriction of $\mathbb{T}_{k/2}$ to $S_{k/2}^{\bot}(N, \chi)$ as a $\ZZ[\zeta_{\varphi(N)}]$-module. In particular the set of operators $T_{p^2}$ for primes $p \leq R$ forms a generating set as an algebra. Moreover, $f \in S_{k/2}(N, \chi)$ is an eigenform for all Hecke operators if and only if it is an eigenform for $T_{p^2}$ for $p \leq R$.

\end{thm}

\section{Proof of the Main Result}

In this section, we will give proof of the main result. Here, main ingredient will be Soma Purkait's result which is a consequence of the Shimura life with the Sturm bound. Note that, we are giving a computerized proof here. 

\begin{proof}
It is well-known that eigenforms will be simultaneous eigenvectors for the associated matrices under Hecke operators. One can implement Hecke operators for half-integral weight modular forms in Magma \cite{Magma} or  Pari/GP \cite{Pari} or SageMath \cite{sage} by using standard definition. For instance, here\footnote{https://tinyurl.com/5n6f7hax} the reader can find the code at Magma for implementing Hecke operators for half-integral weight modular forms.

To prove that our Dedekind-eta quotient is indeed an eigenform, it is suifficient to calculate the value $R$ of Theorem \ref{Soma} and check  whether our Dedekind-eta quotient is an Hecke eigenform for $T_{p^2}$ for $p \leq R$. To do this, one can use the code\footnote{https://tinyurl.com/2s48wfcm} to calculate the value $R$s at Sage version 9.6 and then by the following code there it can be easily seen that all Dedekind eta-quotients listed in Table \ref{Table:1} are indeed Hecke eigenform of half integral weight.

\end{proof}

\end{document}